\documentclass[10pt]{amsart}
\usepackage{amscd}
\usepackage{graphicx}
\usepackage{amsmath, amssymb}
\def\g{\gamma}

\def\b{\beta}
\def\s{\sigma}
\def\a{\alpha}
\def\l{\lambda}
\def\mc{\mathbb{C}}

\newtheorem{theorem}{Theorem}[section]
\newtheorem{proposition}[theorem]{Proposition}
\newtheorem{corollary}[theorem]{Corollary}
\newtheorem{lemma}[theorem]{Lemma}

\begin{document}
\textwidth=13cm

\title{Fundamental group of Desargues configuration spaces}

\thanks{2010 AMS Classification Primary: 20F36, 52C35, 55R80, 57M05; Secondary: 51A20.\\
\indent This research is partially supported by Higher Education Commission, Pakistan.}
\author{ BARBU BERCEANU$^{1}$ ,\,\,SAIMA PARVEEN $^{2}$}
\address{$^{1}$Abdus Salam School of Mathematical Sciences,
 GC University, Lahore-Pakistan, and
 Institute of Mathematics Simion Stoilow, Bucharest-Romania.(Permanent address)}
\email {Barbu.Berceanu@imar.ro}
\address{$^{2}$Abdus Salam School of Mathematical Sciences,
 GC University, Lahore-Pakistan.}
\email {saimashaa@gmail.com}
\maketitle
 \pagestyle{myheadings} \markboth{\centerline {\scriptsize
BARBU BERCEANU,\,\,\,SAIMA PARVEEN   }} {\centerline {\scriptsize
Fundamental group of Desargues configuration spaces }}
\begin{abstract} We compute the fundamental group of various spaces of Desargues configurations in complex
projective spaces: planar and non-planar configurations, with a fixed center and also with an arbitrary center.
\end{abstract}

\maketitle

\section{INTRODUCTION}

Let $M$ be a manifold and $\mathcal{F}_k(M)$ be its {\em ordered
configuration space} of $k$-tuples
$\{(x_1,\ldots,x_k)\in M^k\,|\,x_i\neq x_j,\,\,i\neq j\}$. The $k^{th}$ {\em pure braid group}
of $M$ is the fundamental group of $\mathcal{F}_k(M)$. The pure braid group of the plane, denoted by $%
\mathcal{PB}_n$, has the presentation \cite{F}
$$ \pi_1(\mathcal{F}_n(\mathbb{C}))=\mathcal{PB}_n\cong\langle\alpha_{ij}\,,\,\,1%
\leq i<j\leq n\,\big|\,(YB\,3)_n,(YB\,4)_n\rangle $$
where generators $ \alpha_{ij} $ are represented in the figure and the Yang-Baxter relations
\begin{center}
\begin{picture}(360,60)
\thicklines
\put(70,15){\line(0,1){30}}     \put(67,50){$1$}                 \put(67,0){$1$}
\put(100,15){\line(0,1){30}}    \put(91,50){$i-1$}               \put(91,0){$i-1$}
\put(79,30){$\ldots$}           \put(259,30){$\ldots$}           \put(169,30){$\ldots$}
\put(130,15){\line(0,1){7}}     \put(127,50){$i$}                \put(127,0){$i$}
\put(160,15){\line(0,1){30}}    \put(151,50){$i+1$}              \put(151,0){$i+1$}
\put(190,15){\line(0,1){30}}    \put(183,50){$j-1$}              \put(183,0){$j-1$}
\put(220,45){\line(-6,-1){25}}  \put(185,39.5){\line(-6,-1){20}} \put(130,45){\line(0,-1){16}}
\put(155,34.9){\line(-6,-1){20}}\put(125,31){\line(-6,-1){15}}   \put(30,27){$\a_{ij}$}
\put(110,29){\line(6,-1){45}}   \put(165,21){\line(6,-1){20}}
\put(195,16){\line(6,-1){25}}   \put(217,50){$j$}                \put(217,0){$j$}
\put(250,15){\line(0,1){30}}    \put(240,50){$j+1$}              \put(240,0){$j+1$}
\put(280,15){\line(0,1){30}}    \put(277,50){$n$}                \put(277,0){$n$}
\end{picture}
\end{center}
\noindent $(YB\,3)_n$ and $(YB\,4)_n$ are, for any $1\leq i <j< k\leq n$,
\begin{equation*}
(YB\,3)_n: \a_{ij}\a_{ik}\a_{jk}=\a_{ik}\a_{jk}\a_{ij}=\a_{jk}\a_{ij}\a_{ik}
\end{equation*}
and, for any $1\leq i<j<k<l\leq n$,
$$ (YB\,4)_n: [\a_{kl},\a_{ij}]=[\a_{jl},\a_{jk}^{-1}\a_{ik}\a_{jk}]=[\a_{il},\a_{jk}]=
[\a_{jl},\a_{kl}\a_{ik}\a_{kl}^{-1}]=1. $$
The pure braid group of $S^2\approx\mc {\rm P}^1$ have the
presentation (see \cite{FV} and \cite{F}):
$$\pi_1(\mathcal{F}_{k+1}(S^2))\cong\langle\a_{ij}, 1\leq i<j\leq k\,\big|\,(YB\,3)_k,(YB\,4)_k,D_{k}^2=1\rangle,$$
where $D_{k}=\a_{12}(\a_{13}\a_{23})\ldots(\a_{1k}\ldots\a_{k-1,k})$ (in $\mathcal{B}_k $, the Artin braid group,
$D_k$ is the square of the Garside element $\Delta _k$, see \cite{G} and \cite{B1}). In \cite{B1} we started to study the topology
of configuration spaces under simple geometrical restrictions. Using the geometry of the projective space we can
stratify the configuration space $\mathcal{F}_k(\mc {\rm P}^n)$ with complex submanifolds:
$$\mathcal{F}_k(\mathbb{C}{\rm P}^n)=\mathop{\coprod}\limits_{i=1}^{n} \mathcal{F}_{k}^{i,n}\,\,,$$
where $ \mathcal{F}_{k}^{i,n}$ is the ordered configuration space of all $k$-tuples in $\mc {\rm P}^n$
generating a subspace  of dimension $i$. Their fundamental groups are given by (see \cite{B1}):
\begin{theorem}
\label{th:1} The spaces $\mathcal{F}_{k}^{i,n}$ are simply connected with
the following exceptions

\begin{enumerate}
\item for $k\geq2$,
\begin{equation*}
\pi_1(\mathcal{F}_{k+1}^{1,1})\cong\langle\alpha_{ij},\,\,\,1\leq i<j\leq k\,%
\big|(YB\,3)_k,(YB\,4)_k,D_{k}^2=1\rangle;
\end{equation*}

\item for $k\geq 3$ and $n\geq 2$,
\begin{equation*}
\pi_1(\mathcal{F}_{k+1}^{1,n})\cong\langle\alpha_{ij},\,\,\,1\leq i<j\leq k\,%
\big|(YB\,3)_k,(YB\,4)_k,\,D_{k}=1\rangle.
\end{equation*}
\end{enumerate}
\end{theorem}
In this paper we compute the fundamental groups of various configuration spaces related to projective Desargues
configurations. We do not use special notations for the dual projective space: if $P_1,P_2,P_3$
are three points and $d_1,d_2,d_3$ are three lines in $\mc {\rm P}^2$, $(P_1,P_2,P_3)\in \mathcal{F}_{3}^{1,2}$
is equivalent with the collinearity of these points and $(d_1,d_2,d_3)\in \mathcal{F}_{3}^{1,2}$ is equivalent
with the concurrency of these lines. We define $\mathcal{D}^{2,n}$, the
{\em space of planar Desargues configurations in} $\mc {\rm P}^n$ ($n\geq2$), by
$$\mathcal{D}^{2,n}=\big\{(A_1,B_1,A_2,B_2,A_3,B_3)\in\mathcal{F}_{6}^{2,n}\,\big|\,(d_1,d_2,d_3)\in
\mathcal{F}_{3}^{1,2},\,A_i,B_i\in d_i\setminus\{I\}\big\}$$   
(here $I=d_1\cap d_2\cap d_3$).
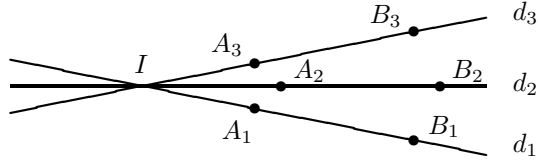
\begin{center}
\begin{picture}(250,60)
\thicklines
\put(20,40){\line(5,-1){180}}   \put(20,20){\line(5,1){180}}      \put(20,30){\line(1,0){180}}
\put(127,33){$A_2$}             \put(187,33){$B_2$}               \put(210,28){$d_2$}
\put(120,27.5){$\bullet$}       \put(180,27.5){$\bullet$}         \put(67,35){$I$}
\put(210,55){$d_3$}             \put(210,5){$d_1$}                \put(110,19){$\bullet$}
\put(100,10){$A_1$}             \put(96,43){$A_3$}                \put(156,54){$B_3$}
\put(178,12){$B_1$}             \put(170,7){$\bullet$}            \put(110,36){$\bullet$}
\put(170,48){$\bullet$}
\end{picture}
\end{center}
We consider also $\mathcal{D}_{I}^{2,n}$, {\em the space of planar Desargues configuration with a
fixed intersection point} $I\in\mc {\rm P}^n$, defined by
$$\mathcal{D}_{I}^{2,n}=\big\{(A_1,B_1,A_2,B_2,A_3,B_3)\in\mathcal{D}^{2,n}\,\big|\,d_1\cap d_2\cap d_3=I\big\}.$$
\begin{theorem}\label{th:1.1} The fundamental group of $\mathcal{D}_{I}^{2,n}$ is given by
$$\pi_1(\mathcal{D}_{I}^{2,n})\cong\left\{
    \begin{array}{ll}
       \mathbb{Z}\oplus\mathbb{Z}\oplus\mathbb{Z}  & \mbox{if }n=2,\\
       \mathbb{Z}\oplus\mathbb{Z}                  & \mbox{if }n\geq3.
    \end{array}
  \right.
$$
\end{theorem}
\noindent The first group is generated by $[\a],[\b]$, and $[\s]$, and the second group is generated by $[\a]$ and $[\b]$.
Precise formulae for $\a,\b$ and $\s$ are given in section $2$; here is a diagram representing these
generators (there is a similar picture for $\b$):
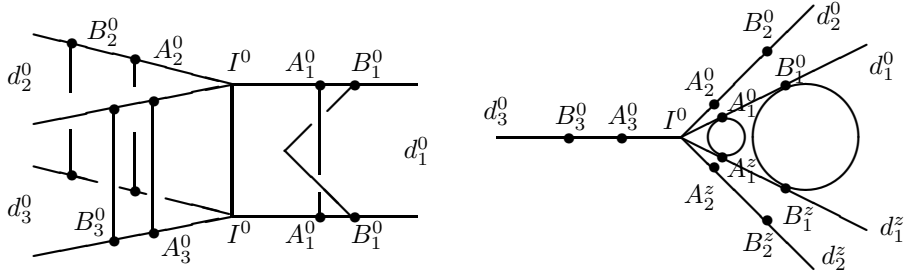
\bigskip
\begin{center}
\begin{picture}(360,100)
\thicklines
\put(90,70){\line(1,0){70}}       \put(90,20){\line(1,0){70}}          \put(110,10){$A_1^0$}
\put(123,20){\line(0,1){9}}       \put(123,70){\line(0,-1){34}}        \put(135,10){$B_1^0$}
\put(110,45){\line(1,-1){25}}     \put(110,45){\line(1,1){10}}         \put(155,42){$d_1^0$}
\put(136,70){\line(-1,-1){10}}    \put(110,75){$A_1^0$}                \put(135,75){$B_1^0$}
\multiput(0,0)(13,0){2}{\multiput(121,17)(0,50){2}{$\bullet$}}         \put(88,75){$I^0$}
\put(90,20){\line(0,1){50}}       \multiput(43,8)(0,50){2}{$\bullet$}  \put(88,10){$I^0$}
\multiput(15,5)(0,50){2}{\line(5,1){75}}
\multiput(45,11)(15,3){2}{\line(0,1){50}}                              \put(35,87){$B_2^0$}
\multiput(29,35)(0,32){2}{\line(0,1){18}}                              \put(60,80){$A_2^0$}
\multiput(58,11)(0,50){2}{$\bullet$}                                   \put(15,89){\line(4,-1){75}}
\multiput(15,39)(50,-12){2}{\line(4,-1){24}}                           \put(30,15){$B_3^0$}
\multiput(0,0)(24,-6){2}{\multiput(27,33)(0,50){2}{$\bullet$}}         \put(63,5){$A_3^0$}
\put(53,29){\line(0,1){23}}       \put(53,79){\line(0,-1){10}}         \put(5,70){$d_2^0$}
\put(47.5,31){\line(4,-1){10}}    \put(5,20){$d_3^0$}                  \put(277,50){\circle{14}}
\put(0,-15){$\a: B_1\mbox{ is moving on the line }d_1^0\setminus\{I^0,A_1^0\}$}
\put(200,-15){$\s: \mbox{ the lines } d_1 \mbox{ and } d_2\mbox{ are moving}$}
\put(190,50){\line(1,0){70}}      \put(253,53){$I^0$}                  \put(212,55){$B_3^0$}
\multiput(215,47)(20,0){2}{$\bullet$}                                  \put(232,55){$A_3^0$}
\put(185,57){$d_3^0$}             \put(260,50){\line(1,1){50}}         \put(260,50){\line(1,-1){50}}
\put(260,50){\line(2,1){70}}      \put(260,50){\line(2,-1){70}}        \put(311,93){$d_2^0$}
\put(313,2){$d_2^z$}              \put(335,13){$d_1^z$}                \put(331,74){$d_1^0$}
\put(261,68){$A_2^0$}             \put(283,90){$B_2^0$}                \put(283,7){$B_2^z$}
\put(277,60){$A_1^0$}             \put(295,73){$B_1^0$}                \put(307,50){\circle{50}}
\multiput(270,60)(20,20){2}{$\bullet$}                                 \put(261,25){$A_2^z$}
\multiput(270,36)(20,-20){2}{$\bullet$}                                \put(298,17){$B_1^z$}
\multiput(273,55)(24,12){2}{$\bullet$}                                 \put(277,35){$A_1^z$}
\multiput(273,40)(24,-12){2}{$\bullet$}
\end{picture}
\end{center}\vspace{1.1cm}
\begin{theorem}\label{th:1}
The fundamental group of $\mathcal{D}^{2,n}$ is given by:
$$\pi_1(\mathcal{D}^{2,n})\cong\left\{
    \begin{array}{ll}
      \mathbb{Z}\oplus\mathbb{Z}    & \mbox{if }n=2,\\
      \mathbb{Z} & \mbox{if }n\geq3.
    \end{array}
  \right.
$$
\end{theorem}
The first group is generated by $[\a]$ and $[\b]$ and the second group is generated by $[\a]$ (or by
$[\b]$); we will use the same notations for $[\a],[\b],[\s]$ and their
images through different natural maps: $ \mathcal{D}_I^{*,*}\to \mathcal{D}^{*,*} $,
$ \mathcal{D}_I^{*,*}\to \mathcal{D}_I^{*,*+1} $, $ \mathcal{D}^{*,*}\to \mathcal{D}^{*,*+1} $.

We define $\mathcal{D}^{3,n}$, the {\em space of non-planar Desargues configurations in} $\mc {\rm P}^n$ ($n\geq3$):
$$\mathcal{D}^{3,n}=\big\{(A_1,B_1,A_2,B_2,A_3,B_3)\in\mathcal{F}_{6}^{3,n}\,\big|\,d_1\cap d_2\cap d_3=I,
\,A_i,B_i\in d_i\setminus\{I\}\big\}$$  
and $\mathcal{D}_{I}^{3,n}$, the associated {\em space of non-planar Desargues configurations
with a fixed intersection point} $I\in\mc {\rm P}^n$.
\begin{theorem}\label{th:1.4}
The fundamental group of $\mathcal{D}_{I}^{3,n}$ is given by:
$$ \pi_1(\mathcal{D}_{I}^{3,n})\cong \left\{
    \begin{array}{ll}
      \mathbb{Z} &   \mbox{if }n=3,    \\
      1          &   \mbox{if }n\geq4 .
    \end{array}
  \right.
$$
\end{theorem}
\begin{theorem}\label{th:1.5}
The fundamental group of $\mathcal{D}^{3,n}$ is given by:
$$\pi_1(\mathcal{D}^{3,n})\cong\left\{
    \begin{array}{ll}
      \mathbb{Z}_4 &   \mbox{if }n=3,    \\
      1            &   \mbox{if }n\geq4.
    \end{array}
  \right.
$$
\end{theorem}
\noindent In the last two theorems, in the non-simply connected cases, the fundamental groups are generated by $[\a]$.

\section{Desargues configurations in the projective plane}

In order to find the fundamental groups of the spaces $\mathcal{D}=\mathcal{D}^{2,2}$ and $\mathcal{D}_{I}=
\mathcal{D}_{I}^{2,2}$ we use two fibrations and their homotopy exact sequences.

\begin{lemma}\label{l.1}
The projection
$$\mu:\mathcal{D}\rightarrow \mc {\rm P}^2,\,(A_1,B_1,A_2,B_2,A_3,B_3)\mapsto I=d_1\cap d_2\cap d_3$$
is a locally trivial fibration with fiber $\mathcal{D}_I.$
\end{lemma}
\begin{proof} Fix a point $I^0\in\mc {\rm P}^2$ and choose a line $l\subset \mc {\rm P}^2\setminus\{I^0\}$ and
the neighborhood $\mathcal{U}_l=\mc {\rm P}^2\setminus l$ of $I^0$. For a point $I$ in this neighborhood and
a Desargues configuration $(A_{1}^0,B_{1}^0,A_{2}^0,B_{2}^0,A_{3}^0,B_3^0)$ on three lines $d_{1}^0,d_{2}^0,d_3^0$
containing $I^0$ construct lines $d_1,d_2,d_3$ containing $I$ and the configuration $(A_1,B_1,\ldots,A_3,B_3)$ as
follows: consider the points $D_i=l\cap d_{i}^0$ and $Q=l\cap I^0I$ and define $d_i=ID_i,\,A_i=d_i\cap QA_{i}^0$ and
in the same way $B_i\,(i=1,2,3)$. We describe this construction using coordinates to show that the map
$$\big(I,(A_{1}^0,B_{1}^0,A_{2}^0,B_{2}^0,A_{3}^0,B_3^0)\big)\mapsto(A_1,B_1,A_2,B_2,A_3,B_3)$$ has a continuous
extension on the singular locus $(d_{1}^0\cup d_{2}^0\cup d_3^0\setminus l)$. Choose a projective frame such
that $I^0=[0:0:1],\, l:X_2=0$. If $I=[s:t:1]$ and $A_{i}^0=[n_i:-m_i:a_i]$, $B_{i}^0=[n_i:-m_i:b_i]$ ($a_i,b_i$ are
distinct and non zero and also $n_im_j\neq m_in_j$ for distinct $i,j=1,2,3$), then we define
$A_{i}=[n_i+sa_i:-m_i+ta_i:a_i]$ and $B_{i}=[n_i+sb_i:-m_i+tb_i:b_i]$, $(i=1,2,3)$, and these formulae agree with the
geometrical construction given for nondegenerate positions of $I\in \mc {\rm P}^2\setminus (d_1^0\cup d_2^0
\cup d_3^0\cup l)$. The trivialization over $\mathcal{U}_l$ is given by
$$\varphi:\mathcal{U}_l\times \mathcal{D}_{I^0}\rightarrow\g^{-1}(\mathcal{U}_{l}),\,\varphi
\big(I,(A_{1}^0,B_{1}^0,A_{2}^0,B_{2}^0,A_{3}^0,B_3^0)\big)=(A_1,B_1,A_2,B_2,A_3,B_3).$$
\end{proof}
\begin{lemma}\label{l:2}
The projection
$$\lambda:\mathcal{D}_{I}\rightarrow\mathcal{F}_{3}(\mc {\rm P}^1),\,(A_1,B_1,A_2,B_2,A_3,B_3)\mapsto(d_1,d_2,d_3)$$
is a locally trivial fibration with fiber $\mathcal{F}_2(\mc )\times\mathcal{F}_{2}(\mc)\times\mathcal{F}_2(\mc)$.
\end{lemma}
\begin{proof}
Fix a point $d_*^0=(d_1^0,d_2^0,d_3^0)$ in $\mathcal{F}_3(\mc {\rm P}^1)$ and choose a point $Q$ in
$\mc {\rm P}^2\setminus(d_{1}^0\cup d_{2}^0\cup d_3^0)$ and the neighborhood $\mathcal{U}_Q=\{(d_1,d_2,d_3)
\in\mathcal{F}_{3}(\mc {\rm P}^1)\big|Q\notin d_1\cup d_2\cup d_3\}$. The trivialization over $\mathcal{U}_Q$ is given by
$$\psi:\mathcal{U}_Q\times \mathcal{F}_{2}(d_{1}^0\setminus\{I\})\times \mathcal{F}_{2}(d_{2}^0
\setminus\{I\})\times\mathcal{F}_2(d_3^0\setminus\{I\})\rightarrow \lambda^{-1}(\mathcal{U}_Q)$$
$$\psi\big((d_1,d_2,d_3),(A_{1}^0,B_{1}^0),(A_{2}^0,B_{2}^0),(A_3^0,B_3^0)\big)=(A_1,B_1,A_2,B_2,A_3,B_3),$$
where $A_i=d_i\cap QA_{i}^0$ and similarly for $B_i\,(i=1,2,3)$. Obviously, $A_i,B_i$ and $I$ are three distinct points on $d_i$.
\end{proof}
In $\mathcal{D}_{I^0=[0:0:1]}$ we choose the base point $D^0=(A_{1}^0,B_{1}^0,A_{2}^0,B_{2}^0,A_{3}^0,B_3^0)$
where, for $k=1,2$, $A_{k}^0=[-1:k:1],B_{k}^0=[-1:k:2], A_{3}^0=[0:1:1],B_3^0=[0:1:2]$.
The corresponding lines are given by the equations $d_k^0:kX_0+X_1=0,\,d_3^0:X_0=0$ and we identify the affine
line $\mathbb{C}$ with $d_k^0$ as follows: for $k=1,2$, $z\mapsto[-1:k:z]$, and for $k=3,\,z\mapsto[0:1:z]$ (therefore
the intersection point $I^0=[0:0:1]$ is the point at infinity of these lines). We  identify the set of three distinct
lines through $I^0$ with the configuration space $\mathcal{F}_{3}(\mc {\rm P}^1)$; in this space the base point
is $d_*^0=(d_1^0,d_2^0,d_3^0)$. In the configuration spaces $\mathcal{F}_2(d_i^0\setminus \{I^0\})$ we choose
the base points $(A_i^0,B_i^0)$, $i=1,2,3$. The homotopy exact sequence from Lemma \ref{l:2} and the triviality
of $\pi_2(\mathcal{F}_3(\mc {\rm P}^1))$ (see \cite{B2}) give the short exact sequence
$$ 1\rightarrow\pi_1(\mathcal{F}_2(\mc))\times \pi_1(\mathcal{F}_2(\mc))\times\pi_1(\mathcal{F}_2(\mc))\mathop
{\rightarrow}\limits^{j_*}\pi_1(\mathcal{D}_{I^0})\mathop{\rightarrow}\limits^{\lambda_*}\pi_1(\mathcal{F}_3(\mc {\rm P}^1))\rightarrow1.$$
{\em Proof of Theorem \ref{th:1.1} (the case $n=2$).}
The first group, isomorphic to  $ \mathbb{Z}^3 $, is generated by the pure braids $a,b,c$, hence their images in
$\pi_1(\mathcal{D}_{I^0})$ are given by the
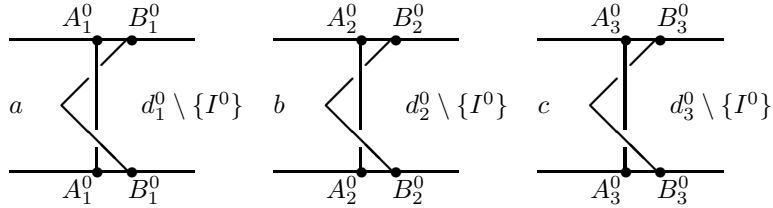
\begin{center}
\begin{picture}(360,80)
\thicklines
\multiput(20,60)(100,0){3}{\line(1,0){70}}      \multiput(20,10)(100,0){3}{\line(1,0){70}}
\multiput(53,10)(100,0){3}{\line(0,1){9}}       \multiput(53,60)(100,0){3}{\line(0,-1){34}}
\multiput(40,35)(100,0){3}{\line(1,-1){25}}     \multiput(40,35)(100,0){3}{\line(1,1){10}}
\multiput(65,60)(100,0){3}{\line(-1,-1){10}}
\multiput(0,0)(100,0){3}{\multiput(0,0)(13,0){2}{\multiput(51,7)(0,50){2}{$\bullet$}}}
\put(70,32){$d_1^0\setminus \{I^0\}$}           \put(40,65){$A_1^0$}                 \put(65,65){$B_1^0$}
\put(20,32){$a$}                                \put(40,0){$A_1^0$}                  \put(65,0){$B_1^0$}
\put(170,32){$d_2^0\setminus \{I^0\}$}          \put(140,65){$A_2^0$}                \put(165,65){$B_2^0$}
\put(120,32){$b$}                               \put(140,0){$A_2^0$}                 \put(165,0){$B_2^0$}
\put(270,32){$d_3^0\setminus \{I^0\}$}          \put(240,65){$A_3^0$}                \put(265,65){$B_3^0$}
\put(220,32){$c$}                               \put(240,0){$A_3^0$}                 \put(265,0){$B_3^0$}
\end{picture}
\end{center}
\bigskip
\noindent homotopy classes of the maps $\a,\b,\g: (S^1,1)\rightarrow (\mathcal{\mathcal{D}}_{I^0},D^0)$
$$\begin{array}{ll}
   \a(z)=(A_{1}^{0},B_{1}^{\a(z)},A_2^0,B_2^0,A_3^0,B_3^0),    & B_1^{\a(z)}=[-1:1:1+z],\\
   \b(z)=(A_{1}^{0},B_{1}^{0},A_2^0,B_2^{\b(z)},A_3^0,B_3^0),  & B_2^{\b(z)}=[-1:2:1+z],\\
   \g(z)=(A_{1}^{0},B_{1}^{0},A_2^0,B_2^0,A_3^0,B_3^{\g(z)}),  & B_3^{\g(z)}=[0:1:1+z].
\end{array}$$
The third group, $\pi_1(\mathcal{F}_3(\mc {\rm P}^1)\cong \mathbb{Z}_2$, is generated by the homotopy class of the map
$$s: (S^1,1)\rightarrow (\mathcal{F}_3(\mc {\rm P}^1),d_*^0),\,z\mapsto (d_1^{s(z)}:z X_0+X_1=0, d_2^{s(z)}:2z X_0+X_1=0, d_3^0 ),$$
because this corresponds to the braid $\a_{12}$ in $\mc {\rm P}^1$. We lift the map $s$ to the map
$$\s: (S^1,1)\rightarrow (\mathcal{D}_I^0,D^0),\,z\mapsto(A_{1}^{\s(z)},B_{1}^{\s(z)},A_2^{\s(z)},B_2^{\s(z)},A_3^0,B_3^0),$$
where $A_{k}^{\s(z)}=[-1:kz:1]$, $B_{k}^{\s(z)}=[-1:kz:2]$, $k=1,2$.

The group $\pi_1(\mathcal{D}_{I^0},D^0)$ is generated by the homotopy classes of $\a,\b,\g$ and $\s$; the
defining relations are commutation relations between $[\a],[\b]$ and $[\g]$ from $\pi_1(\mathcal{F}_2(\mathbb{C})^3)$
and the four relations, to be proved in the next two lemmas:
$$ \begin{array}{ll}
   \a)    &  [\s][\a][\s]^{-1}=[\a],\\
   \b)    &  [\s][\b][\s]^{-1}=[\b],\\
   \g)    &  [\s][\g][\s]^{-1}=[\g],\\
   \s)    &  [\s]^2=[\a]^{-1}[\b]^{-1}[\g].
\end{array}
$$
The generator $[\g]$ can be eliminated, $[\s]$ commutes with $[\a]$ and $[\b]$, and the third relation,
$\g)$, is a consequence of the previous commutation relations. \hfill $\square$
\begin{lemma}
In $\pi_1(\mathcal{D}_{I^0},D^0)$ the next relation holds:

$\s)\,\,[\s]^2=[\a]^{-1}[\b]^{-1}[\g].$
\end{lemma}
\begin{proof} The map
$$\Lambda:(D^2,S^1)\rightarrow(\mathcal{F}_3(\mc {\rm P}^1),d_*^0=(d_1^0,d_2^0,d_3^0)),\,z\mapsto(d_1^{\Lambda(z)},d_2^{\Lambda(z)},d_3^{\Lambda(z)}),$$
where $d_k^{\Lambda(z)}:(kz-r)X_0+(\overline{z}+kr)X_1=0,\,(k=1,2)$, and $d_3^{\Lambda(z)}: zX_0+rX_1=0$ (the notation
 $r=1-|z|$ will be used in this proof and the next ones),
shows that $s^2\simeq\hbox{constant}_{d_{*}^0}$. We lift this homotopy to
$$\widetilde{\Lambda}: D^2\rightarrow \mathcal{D}_{I^0},\,\widetilde{\Lambda}(z)=\big(A_1^{\widetilde{\Lambda}(z)},B_1^{\widetilde{\Lambda}(z)},
A_2^{\widetilde{\Lambda}(z)},B_2^{\widetilde{\Lambda}(z)},A_3^{\widetilde{\Lambda}(z)},B_3^{\widetilde{\Lambda}(z)}\big),$$
where $A_k^{\widetilde{\Lambda}(z)}=[-\overline{z}-kr:kz-r:\overline{z}],\,B_k^{\widetilde{\Lambda}(z)}=
[-\overline{z}-kr:kz-r:\overline{z}+1],\,(k=1,2)$, and
$A_3^{\widetilde{\Lambda}(z)}=[-r:z:z],B_3^{\widetilde{\Lambda}(z)}=[-r:z:z+1]$; the map
$$ \widetilde{\Lambda}|_{S^1}: S^1\rightarrow\mathcal{D}_{I^0},\,z\mapsto(A_1^z,B_1^z,A_2^z,B_2^z,A_3^0,B_3^z)$$
(with $A_k^z=[-1:kz^2:1]$, $B_k^z=[-1:kz^2:1+z]$, $k=1,2$, and $B_3^z=[0:1:1+\overline{z}]$) has a
trivial homotopy class, therefore we have the relation $[\s]^2=[\s*\s*(\widetilde{\Lambda}|_{S^1})^{-1}]$.

Now we construct a homotopy between $\s*\s*(\widetilde{\Lambda}|_{S^1})^{-1}$ and $\a^{-1}*\b^{-1}*\g$:
$$L: S^1\times I\rightarrow\mathcal{D}_{I^0},\,(z,t)\mapsto\big(A_1^{L(z,t)},B_1^{L(z,t)},A_2^{L(z,t)},B_2^{L(z,t)},
A_3^0,B_3^{L(z,t)}\big),$$
where ($k=1,2$):
$$ A_k^{L(z,t)}=[-1:kL^1(z,t):1], B_k^{L(z,t)}=[-1:kL^1(z,t):L_k^2(z,t)] \, $$
$$ B_3^{L(z,t)}= \left\{   \begin{array}{l}
                                   \,[0:1:2]       \\
                                   \,[0:1:1+z^2]
                        \end{array}
              \right.
                        \begin{array}{l}
                               0\leq \mbox{arg}\,z\leq \pi\\
                               \pi \leq \mbox{arg}\,z\leq 2\pi,
                        \end{array}     $$
and
$$ \begin{array}{llll}
  L^1(z,t) & = & \left\{   \begin{array}{l}
                                    \,z^4           \\
                                    \,\exp(4t\pi i) \\
                                    \,\overline{z}^4
                              \end{array}
                    \right.
                            & \begin{array}{l}
                               0\leq \mbox{arg}\,z\leq t\pi                \\
                               t\pi \leq \mbox{arg}\,z\leq (2-t)\pi        \\
                               (2-t)\pi \leq \mbox{arg}\,z\leq 2\pi,
                              \end{array}                                  \\
  L_k^2(z,t) & = & \left\{    \begin{array}{l}
                                    \,2            \\
                                    \,1+\exp\big(4\dfrac{(2-k)t\pi-\mbox{arg}\,z}{1+t}i\big)  \\
                                    \,2
                              \end{array}
                    \right.
                            & \begin{array}{l}
                              0\leq\mbox{arg}\,z\leq \frac{t+k-1}{k}\pi \\
                             \frac{t+k-1}{k}\pi\leq\mbox{arg}\,z\leq \frac{1+(5-2k)t}{3-k}\pi  \\
                             \frac{1+(5-2k)t}{3-k}\pi\leq \mbox{arg}\,z\leq 2\pi.
                             \end{array}
  \end{array}$$
It is easy to check that $L(-,0)=(\a^{-1}*\b^{-1})*\g$ and $L(-,1)=(\s*\s)*(\widetilde{\Lambda}|_{S^1})^{-1}$.
\end{proof}
\begin{lemma}
In $\pi_1(\mathcal{D}_{I^0},D^0)$ the next relations hold:

$ \a)\,\, [\s][\a][\s]^{-1}=[\a];$

$ \b)\,\, [\s][\b][\s]^{-1}=[\b];$

$ \g)\,\, [\s][\g][\s]^{-1}=[\g].$
\end{lemma}
\begin{proof} The loop $\s*\a*\s^{-1}$ in $\mathcal{D}_{I^0}$ is given by
$z\mapsto (A_1^{\widetilde{\a}(z)},B_1^{\widetilde{\a}(z)},A_2^{\widetilde{\a}(z)},B_2^{\widetilde{\a}(z)},A_3^0,B_3^0)$,
where the points $A_k^{\widetilde{\a}(z)}$ ($k=1,2$), $B_1^{\widetilde{\a}(z)}$ and $B_2^{\widetilde{\a}(z)}$ are given by :
$$\begin{array}{llll}
A_k=[-1:kz^3:1]            & B_1=[-1:z^3:2]            & B_2=[-1:2z^3:2]
                                                       &  \mbox{arg}\,z\in[0,\frac{2\pi}{3}]              \\
A_k=A_{k}^0                & B_1=[-1:1:1+z^3]          & B_2=B_2^0
                                                       & \mbox{arg}\,z\in[\frac{2\pi}{3},\frac{4\pi}{3}] \\
A_k=[-1:k\overline{z}^3:1] & B_1=[-1:\overline{z}^3:2] & B_2=[-1:2\overline{z}^3:2]
                                                       & \mbox{arg}\,z\in [\frac{4\pi}{3},2\pi].
\end{array}   $$
We define two maps
$$ \varepsilon: S^1\times I\rightarrow S^1\,,\,\,\varepsilon(z,t)=\left\{
     \begin{array}{ll}
  z^3            & 0\leq\mbox{arg}\,z\leq\frac{2t}{3}\pi \\
  \exp(2t\pi i)   & \frac{2t}{3}\pi \leq \mbox{arg}\,z\leq \frac{2(3-t)}{3}\pi\\
  \overline{z}^3 & \frac{2(3-t)}{3}\pi\leq\mbox{arg}\,z\leq 2\pi,
     \end{array}
  \right.  $$

$$ \eta:S^1\rightarrow \mc \setminus\{1\}\,,\, \eta(z)=\left\{
     \begin{array}{ll}
  2     & \mbox{arg}\,z\in[0,\frac{2\pi}{3}]\cup[\frac{4\pi}{3},2\pi]\\
  1+z^3 & \mbox{arg}\,z\in[\frac{2\pi}{3},\frac{4\pi}{3}].
      \end{array}
  \right.  $$
and a new homotopy
$$ K_{\a}(z,t): S^1\times I\rightarrow \mathcal{D}_{I^0},\,K_{\a}(z,t)=\big(A_1(z,t),
\widetilde{B}_1(z,t),A_2(z,t),B_2(z,t),A_3^0,B_3^0\big),$$
where $A_k(z,t)=[-1:k\varepsilon(z,t):1]$, $B_k(z,t)=[-1:k\varepsilon(z,t):2]$, ($k=1,2$),
$\widetilde{B}_1(z,t)=[-1:\varepsilon(z,t):\eta(z)]$. One can check that $K_\a|_{t=0}\simeq \a $
and $K_\a|_{t=1}=\s*\a*\s^{-1}$. Similarly we have a homotopy $K_\b$
between $\b$ and $K_\b|_{t=1}=\s*\b*\s^{-1}$. Next homotopy (we also use the notation $ B_3(z,t)=[0:1:\eta(z)]$)
$$ K_{\g}(z,t): S^1\times I\rightarrow \mathcal{D}_{I^0},\,(z,t)\mapsto\big(A_1(z,t),   
B_1(z,t),A_2(z,t),B_2(z,t),A_3^0,B_3(z,t)\big), $$
gives the last relation: $ K_{\g|t=0}\simeq \g $, $ K_{\g|t=1}=\s*\g*\s^{-1} $.
\end{proof}

\noindent{\em Proof of Theorem \ref{th:1} (the case $n=2$)}. Lemma \ref{l.1} gives the exact sequence
$$\ldots\longrightarrow \pi_2(\mc {\rm P}^2)\mathop{\longrightarrow}\limits^{\delta_*} \mathbb{Z}\oplus
\mathbb{Z}\oplus\mathbb{Z}\longrightarrow \pi_1(\mathcal{D})\longrightarrow 1$$
where the first group is cyclic generated by the homotopy class of the map
$$\Phi:(D^2,S^1)\rightarrow(\mc {\rm P}^2,I^0),\,z\mapsto[0:r:z].$$
We choose the lift
$$ \widetilde{\Phi}:(D^2,S^1)\rightarrow(\mathcal{D},\mathcal{D}_{I^0}),\,z\mapsto\big(A_1
^{\widetilde{\Phi}(z)},B_1^{\widetilde{\Phi}(z)},A_2^{\widetilde{\Phi}(z)},
B_2^{\widetilde{\Phi}(z)},A_3^{\widetilde{\Phi}(z)},B_3^{\widetilde{\Phi}(z)}\big), $$
where ($k=1,2$)
$$ \begin{array}{lll}
A_k^{\widetilde{\Phi}(z)}  & = & \big[-1:(2k+1)r+k\overline{z}:(2k+1)z+k(r-2)\big],\\
B_k^{\widetilde{\Phi}(z)}  & = & \big[-1:(2k+2)r+k\overline{z}:(2k+2)z+k(r-2)\big],\\
A_3^{\widetilde{\Phi}(z)}  & = & \big[-r:\overline{z}+4r:4z-3(r+1)\big],\\
B_3^{\widetilde{\Phi}(z)}  & = & \big[-r:\overline{z}+5r:5z-3(r+1)\big],
\end{array}  $$
hence ${\rm Im}\, \delta_*$ is generated by the homotopy class of the map
$$ \widetilde{\Phi}\big|_{S^1}:S^1\rightarrow \mathcal{D}_{I^0}\,,\,z\mapsto\big(A_1^{\Phi(z)},B_1^{\Phi(z)},
A_2^{\Phi(z)},B_2^{\Phi(z)},A_3^{\Phi(z)},B_3^{\Phi(z)}\big),$$
with $(k=1,2)$
$$ \begin{array}{ll}
A_k^{\Phi(z)}=[-1:k\overline{z}:(2k+1)z-2k], & B_k^{\Phi(z)}=[-1:k\overline{z}:(2k+2)z-2k], \\
A_3^{\Phi(z)}=[0:\overline{z}:4z-3],         & B_3^{\Phi(z)}=[0:\overline{z}:5z-3].
\end{array}$$
The maps $\l\circ\widetilde{\Phi}\big|_{S^1}$ and $s^{-1}$ coincide, therefore the product
$[\widetilde{\Phi}\big|_{S^1}]\cdot[\s]$ belongs to ${\rm ker}\, \lambda_*={\rm Im}\, j_*$. We show
that $[\widetilde{\Phi}\big|_{S^1}]\cdot[\s]=[\a]\cdot[\b]\cdot[\g]$ and this implies the claim
of the theorem. We define the homotopy:
$$ H:S^1\times I\rightarrow \mathcal{D}_{I^0},\, (z,t)\mapsto \big(A_1^{H(z,t)},B_1^{H(z,t)},
A_2^{H(z,t)},B_2^{H(z,t)},A_3^{H(z,t)},B_3^{H(z,t)}\big),$$
where $(k=1,2)$
$$ \begin{array}{ll}
      A_k^{H(z,t)} =[-1:H_k^1(z,t):H_k^2(z,t)] &  B_k^{H(z,t)}=[-1:H_k^1(z,t):H_k^2(z,t)+H_k^4(z,t)]  \\
      A_3^{H(z,t)} =[0:1:H^3(z,t)] &  B_3^{H(z,t)}=[0:1:H^3(z,t)+H^5(z,t)]
   \end{array}  $$
and
$$ \begin{array}{llll}
       H_k^1(z,t) & = & \left\{  \begin{array}{l}
                                    k\overline{z}^2    \\
                                    k\exp(-2t\pi i)     \\
                                    kz^2
                                  \end{array}
                        \right. & \begin{array}{l}
                                    0\leq \mbox{arg}\,z\leq t\pi        \\
                                    t\pi\leq \mbox{arg}\,z\leq (2-t)\pi \\
                                    (2-t)\pi\leq \mbox{arg}\,z\leq 2\pi,
                                   \end{array}   \\
       H_k^2(z,t) & =  & \left\{  \begin{array}{l}
                                    1+(2k+1)t(z^2-1)  \\
                                    1
                                  \end{array}
                         \right.
                                &  \begin{array}{l}
                                      0\leq \mbox{arg}\,z\leq \pi        \\
                                      \pi\leq \mbox{arg}\,z\leq 2\pi,
                                   \end{array}    \\
       H^3(z,t)  &  =  & \left\{  \begin{array}{l}
                                      1+t(4z^4-3z^2-1)     \\
                                      1
                                   \end{array}
                         \right.
                                &  \begin{array}{l}
                                      0\leq \mbox{arg}\,z\leq \pi\\
                                      \pi\leq \mbox{arg}\,z\leq 2\pi,
                                   \end{array}   \\
       H_1^4(z,t)  &  =  & \left\{ \begin{array}{l}
                                      \exp\big(\dfrac{4\mbox{arg}\,z}{1+t}i\big)                        \\
                                      1
                                    \end{array}
                         \right.
                                &  \begin{array}{l}
                                      0\leq \mbox{arg}\,z\leq \frac{1+t}{2}\pi\\
                                      \frac{1+t}{2}\pi\leq \mbox{arg}\,z\leq 2\pi,
                                   \end{array} \\
       H_2^4(z,t)  &  =  & \left\{ \begin{array}{l}
                                      1                           \\
                                      \exp\big(2\dfrac{2\mbox{arg}\,z-(1-t)\pi}{1+t}i\big) \\
                                      1
                                     \end{array}
                            \right.
                                 &   \begin{array}{l}
                                       0\leq \mbox{arg}\,z\leq \frac{1-t}{2}\pi         \\
                                       \frac{1-t}{2}\pi\leq \mbox{arg}\,z\leq \pi  \\
                                       \pi\leq \mbox{arg}\,z\leq 2\pi .
                                   \end{array}  \\
       H^5(z,t)  &  =  & \left\{ \begin{array}{l}
                                      1                            \\
                                      \exp[4(\mbox{arg}\,z-(1-t)\pi)i] \\
                                      1
                                     \end{array}
                            \right.
                                 &   \begin{array}{l}
                                       0\leq \mbox{arg}\,z\leq (1-t)\pi         \\
                                       (1-t)\pi\leq \mbox{arg}\,z\leq (2-t)\pi  \\
                                       (2-t)\pi\leq \mbox{arg}\,z\leq 2\pi .
                                   \end{array}  \\
\end{array}  $$
These computations give ${\rm Im}\, \delta_*=\mathbb{Z}\langle\,2[\a]+2[\b]+[\s]\,\rangle$, therefore we can
choose $[\a]$ and $[\b]$ as generators of the fundamental group of $\mathcal{D} $. \hfill $\square$

\section{Planar Desargues configuration in $\mc {\rm P}^n$}

First we reduce the computations of $\pi_1(\mathcal{D}_{I}^{2,n})$ and of $\pi_1(\mathcal{D}^{2,n})$ to the case $n=3$.
\begin{lemma}\label{lem.3.1}
The following projections are locally trivial fibrations:
$$  a)\,\,\,\,\mathcal{D}_{I}^{2,2}\hookrightarrow\mathcal{D}_{I}^{2,n}\rightarrow {\rm Gr}^1(\mc {\rm P}^{n-1}),\,\,\,
(A_1,B_1,A_2,B_2,A_3,B_3)\mapsto \hbox{line}\,(d_1,d_2,d_3)\, ;   $$
$$ b)\,\,\,\mathcal{D}^{2,2}\hookrightarrow\mathcal{D}^{2,n}\rightarrow {\rm Gr}^2(\mc {\rm P}^{n}),\,\,\,
(A_1,B_1,A_2,B_2,A_3,B_3)\mapsto \hbox{2-plane}\,(d_1,d_2,d_3). $$
\end{lemma}
\begin{proof} a) Fix a 2-plane $P_0$ through $I$ and choose a hyperplane $H\subset \mc {\rm P}^n$ such
that $I\notin H $ and an $(n-3)$ dimensional subspace $Q\subset H$ such that $Q\cap l_0=\emptyset$,
where $l_0=P_0\cap H$. Take as a neighborhood of $P_0$ the set $\{P\hbox{ a 2-plane in }\,\mc {\rm P}^n\mid I\in P,
P\cap Q=\emptyset \}$ and associate to a Desargues configuration in $\mathcal{D}_{I}(P_0)$ the projection from $Q$,
an element in $\mathcal{D}_{I}(P)$: $C_i^0=d_{i}^0\cap l_0$, $l=P\cap H$, $C_i=(Q\vee C_i^0)\cap l$,
$Q_i=Q\cap(C_i C_i^0)$, $d_i=IC_i $, $A_i=Q_iA_i^0\cap d_i$, $B_i=Q_iB_i^0\cap d_i$ (for $i=1,2,3$). Using projective coordinates
one can show that this trivialization is well defined on the singular locus $P=P_0$: if $ I=[0:\ldots:0:1]$,
$P_0:X_0=\ldots=X_{n-3}=0$, $A_i^0=[0:\ldots:a_{n-2,i}^0:a_{n-1,i}^0:a_{n,i}^0] $, $B_i^0=[0:\ldots:b_{n-2,i}^0:b_{n-1,i}^0:b_{n,i}^0] $,
and $P$ is defined by the equations $X_k=p_{k,1}X_{n-2}+p_{k,2}X_{n-1}+p_{k,3}X_n$ $(k=0,\ldots,n-3)$, then
$A_i=[p_{0,0}a_{n-2,i}+p_{0,1}a_{n-1,i}:\ldots:p_{n-3,0}a_{n-2,i}+p_{n-3,1}a_{n-1,i}:a_{n-2,i}^0:a_{n-1,i}^0:a_{n,i}^0] $, $B_i=[p_{0,0}b_{n-2,i}+p_{0,1}b_{n-1,i}:\ldots:p_{n-3,0}b_{n-2,i}+p_{n-3,1}b_{n-1,i}:b_{n-2,i}^0:b_{n-1,i}^0:b_{n,i}^0] $.

b) Fix a 2-plane $P_0$ and choose as center of projection a disjoint $n-3$ dimensional subspace $Q$.
Take as a neighborhood of $P_0$ the set of 2-planes disjoint from $Q$. The projection from $Q$ associate to a
Desargues configuration in $\mathcal{D}^2(P_0)$ a Desargues configuration in
$\mathcal{D}^2(P): P\cap(Q\vee I^0)=I,\,P\cap(Q\vee d_i^0)=d_i,\,d_i\cap (Q\vee A_i^0)=A_i,\,d_i\cap(Q\vee B_i^0)=B_i$.
\end{proof}
\begin{corollary}\label{cor:3.2}For $n\geq3$ we have

$ a)\,\,\,\pi_1(\mathcal{D}_{I}^{2,3})\cong \pi_1(\mathcal{D}_{I}^{2,n});$

$ b)\,\,\,\pi_1(\mathcal{D}^{2,3})\cong \pi_1(\mathcal{D}^{2,n}).$
\end{corollary}
\begin{proof} This is a consequence of the stability of the second homotopy group of the complex Grassmannians:
\begin{center}
\begin{picture}(360,65)
\thicklines
\put(40,50){$\pi_2({\rm Gr}^1(\mc {\rm P}^2))$}  \put(110,53){\vector(1,0){20}}   \put(135,50){$\pi_1(\mathcal{D}_I^{2,2})$}
\put(179,53){\vector(1,0){20}}       \put(203,50){$\pi_1(\mathcal{D}_I^{2,3})$}   \put(245,53){\vector(1,0){20}}
\put(270,48){$1$}                    \put(220,41){\vector(0,-1){15}}              \put(110,13){\vector(1,0){20}}
\put(135,10){$\pi_1(\mathcal{D}_I^{2,2})$} \put(179,13){\vector(1,0){20}}         \put(203,10){$\pi_1(\mathcal{D}_I^{2,n})$}
\put(245,13){\vector(1,0){20}}       \put(270,8){$1$}                             \put(79,41){\vector(0,-1){15}}
\put(85,30){$\cong$}                 \put(150,41){\vector(0,-1){15}}              \put(155,30){$\cong$}
\put(30,10){$\pi_2({\rm Gr}^1(\mc {\rm P}^{n-1}))$}
\end{picture}
\end{center}
and also
\begin{center}
\begin{picture}(360,65)
\thicklines
\put(40,50){$\pi_2({\rm Gr}^2(\mc {\rm P}^3))$}   \put(110,53){\vector(1,0){20}}       \put(135,50){$\pi_1(\mathcal{D}^{2,2})$}
\put(179,53){\vector(1,0){20}}        \put(203,50){$\pi_1(\mathcal{D}^{2,3})$}         \put(245,53){\vector(1,0){20}}
\put(270,48){$1$}                     \put(40,10){$\pi_2({\rm Gr}^2(\mc {\rm P}^n))$}  \put(110,13){\vector(1,0){20}}
\put(135,10){$\pi_1(\mathcal{D}^{2,2})$}  \put(179,13){\vector(1,0){20}}               \put(203,10){$\pi_1(\mathcal{D}^{2,n})$}
\put(245,13){\vector(1,0){20}}        \put(270,8){$1\, .$}                             \put(79,41){\vector(0,-1){15}}
\put(85,30){$\cong$}                  \put(150,41){\vector(0,-1){15}}                  \put(155,30){$\cong$}
\put(220,41){\vector(0,-1){15}}
\end{picture}
\end{center}
\end{proof}
Using the fibration of Lemma \ref{lem.3.1}$\,a)$ for $n=3$ we have the exact sequence
$$\ldots\rightarrow\pi_2(\mc {\rm P}^2)\mathop{\rightarrow}\limits^{\delta_*}\pi_1(\mathcal{D}_I^{2,2})
\rightarrow\pi_1(\mathcal{D}_I^{2,3})\rightarrow1.$$
We choose the base point in $\mathcal{D}_I^{2,3}$ the image of the base point in $\mathcal{D}_I$ through
the embedding $[x_0:x_1:x_2]\mapsto[x_0:x_1:x_2:0]$ and we denote the compositions
$\a,\b: S^1\rightarrow\mathcal{D}_I^{2,2}\rightarrow\mathcal{D}_I^{2,3}$ with the same letters.
\begin{proposition}\label{pr.3.1} In the exact sequence of the fibration $\mathcal{D}_I^{2,3}\rightarrow \mc {\rm P}^{2}$ we have:

a) ${\rm Im}\, \delta_*=\mathbb{Z}([\a]+[\b]+[\s])$;

b) $\pi_1(\mathcal{D}_I^{2,3})\cong\mathbb{Z}\oplus\mathbb{Z}$ is generated by $[\a]$ and $[\b]$.
\end{proposition}
\begin{proof} a) The base point in ${\rm Gr}^1(\mc {\rm P}^2)\approx \mc {\rm P}^2$ is the line $X_3=0$ (in the dual space
of lines through $I^0=[0:0:1:0]$) and we choose the generator of $\pi_2(\mc {\rm P}^2)$ the homotopy class of the map
$$\Pi:(D^2,S^1)\rightarrow {\rm Gr}^1(\mc {\rm P}^2),\,z\mapsto (1-|z|)X_1+zX_3=0.$$
The lift $\widetilde{\Pi}:D^2\rightarrow \mathcal{D}_{I^0}^{2,3},\,z\mapsto \big(A_1^{\widetilde{\Pi}(z)},
B_1^{\widetilde{\Pi}(z)},A_2^{\widetilde{\Pi}(z)},B_2^{\widetilde{\Pi}(z)},
A_3^{\widetilde{\Pi}(z)},B_3^{\widetilde{\Pi}(z)}\big)$ is given by ($k=1,2$)
$$\begin{array}{ll}
A_k^{\widetilde{\Pi}(z)}=[2r|z|-1:kz:1:-kr],  & A_3^{\widetilde{\Pi}(z)}=[0:z:z:-r], \\
B_k^{\widetilde{\Pi}(z)}=[2r|z|-1:kz:2:-kr],  & B_3^{\widetilde{\Pi}(z)}=[0:z:z+1:-r],
\end{array} $$
where the corresponding lines are
$$d_k^{\widetilde{\Pi}(z)}:kX_0+\overline{z}X_1-rX_3=0,\,\,rX_1+zX_3=0, \,d_3^{\widetilde{\Pi}(z)}: X_0=0,\,rX_1+zX_3=0.$$
The homotopy
$$ M:S^1\times I\rightarrow \mathcal{D}_{I^0}^{2,2}, \, (z,t)\mapsto \big( A_1^{M(z,t)},B_1^{M(z,t)},
A_2^{M(z,t)},B_2^{M(z,t)},A_3^0,B_3^{M(z,t)}\big),$$
where $ A_k^{M(z,t)}=[-1:km_1(z,t):1]$, $B_k^{M(z,t)}=[-1:km_1(z,t):2]$, and $B_3^{M(z,t)}=[0:1:1+m_2(z,t)] $ are defined by:
$$ \begin{array}{llll}
       m_1(z,t) & = & \left\{  \begin{array}{l}
                                    \exp\big(2\dfrac{\mbox{arg}\, z}{2-t}i\big)    \\
                                    1
                                  \end{array}
                        \right. & \begin{array}{l}
                                    0\leq \mbox{arg}\,z\leq (2-t)\pi        \\
                                    (2-t)\pi\leq \mbox{arg}\,z\leq 2\pi,
                                   \end{array}   \\
       m_2(z,t) & =  & \left\{  \begin{array}{l}
                                    1  \\
                                    exp\big(2\dfrac{t\pi-\mbox{arg}\, z}{2-t}i\big)
                                  \end{array}
                         \right.
                                &  \begin{array}{l}
                                      0\leq \mbox{arg}\,z\leq t\pi        \\
                                      t\pi\leq \mbox{arg}\,z\leq 2\pi,
                                   \end{array}
  \end{array} $$
shows that the restriction $\widetilde{\Pi}|_{S^1}$ and the loop $\s*\g^{-1}$ are homotopic. Using this and
the relation $[\g]=[\a]+[\b]+2[\s]$ we find $\delta_*([\Pi])=[\widetilde{\Pi}|_{S^1}]=-[\a]-[\b]-[\s]$.

b) The second part is a consequence of part a).
\end{proof}
\begin{proposition}\label{pr.3.2}
The fundamental group of $\mathcal{D}^{2,3}$ is isomorphic to $\mathbb{Z} $ and it is generated
by $[\a]$ (or by $[\b]$).
\end{proposition}
\begin{proof} This is a consequence of Proposition \ref{pr.3.1} and the computations in section 2:
\begin{center}
\begin{picture}(360,100)
\thicklines
\put(10,85){$\pi_2(\mc {\rm P}^2)=\mathbb{Z}\langle[\Phi]\rangle$}  \put(100,88){\vector(1,0){20}}   \put(105,92){$\delta_*^2$}
\put(130,85){$\pi_1(\mathcal{D}_I^{2,2})=\mathbb{Z\langle[\a],[\b],[\s]\rangle}$}                \put(245,88){\vector(1,0){20}}
\put(270,85){$\pi_1(\mathcal{D}^{2,2})$}           \put(315,88){\vector(1,0){20}}  \put(340,83){$1$}
\put(10,45){$\pi_2(\mc {\rm P}^3)=\mathbb{Z}\langle[\Phi^3]\rangle$}  \put(100,48){\vector(1,0){20}}   \put(105,52){$\delta_*^3$}
\put(130,45){$\pi_1(\mathcal{D}_I^{2,3})=\mathbb{Z\langle[\a],[\b]\rangle}$}                     \put(245,48){\vector(1,0){20}}
\put(270,45){$\pi_1(\mathcal{D}^{2,3})$}         \put(315,48){\vector(1,0){20}}    \put(340,43){$1$}
\put(50,76){\vector(0,-1){15}}                   \put(55,65){$\cong$}              \put(180,76){\vector(0,-1){15}}
\put(185,65){$i_*$}                              \put(290,76){\vector(0,-1){15}}   \put(295,65){$i_*$}
\end{picture}
\end{center}
\vspace{-1cm}
hence $\delta_*^3([\Phi^3])=i_*\delta_*^2([\Phi])=i_*([\widetilde{\Phi}_{|S^1}])=i_*(2[\a]+2[\b]+[\s])=[\a]+[\b]$.
\end{proof}

\section{Non planar Desargues Configurations}

First we analyze the fundamental group of two three-dimensional configuration spaces
$\mathcal{D}_I^3=\mathcal{D}_I^{3,3}$ and $\mathcal{D}^3=\mathcal{D}^{3,3}$.
\begin{lemma}\label{lem.4.1}
The following projections are locally trivial fibrations:\\
$a)\,\,\mathcal{F}_2(\mc)\times \mathcal{F}_2(\mc)\times \mathcal{F}_2(\mc)\hookrightarrow\mathcal{D}_I^3\rightarrow\mathcal{F}_3^{2,2},\,(A_1,B_1,A_2,B_2,A_3,B_3)\mapsto(d_1,d_2,d_3)$\\
$b)\,\,\mathcal{D}_I^3\hookrightarrow\mathcal{D}^3\rightarrow\mc {\rm P}^3,\,(A_1,B_1,A_2,B_2,A_3,B_3)\mapsto I=d_1\cap d_2\cap d_3$.
\end{lemma}
\begin{proof} The proofs are similar to those of Lemmas 2.1 and 2.2.
\end{proof}
\noindent{\em Proof of Theorem \ref{th:1.4} (the case $n=3$) }. We modify a little the previous notations:
the base point in these solid Desargues configurations are related to the center $I^0=[0:0:1:0]$ and to the points:
$$A_1^0=[0:0:0:1],\, B_1^0=[0:0:1:1],\, d_1^0:X_0=X_1=0,$$
$$A_2^0=[0:1:0:0],\, B_2^0=[0:1:1:0],\, d_2^0:X_0=X_3=0,$$
$$A_3^0=[1:0:0:0],\, B_3^0=[1:0:1:0],\, d_1^0:X_1=X_3=0.$$
\noindent Using the fibrations of Lemma \ref{lem.4.1} we find
$$ \pi_2(\mathcal{F}_3^{2,2})\mathop{\rightarrow}\limits^{\delta_*}\pi_1(\mathcal{F}_2(\mc)^3)\cong\mathbb{Z}^3
\rightarrow\pi_1(\mathcal{D}_{I^0}^3)\rightarrow1, $$
where the first group is isomorphic with $\pi_2\big(\mathcal{F}_2(\mc {\rm P}^2)\big)\cong \mathbb{Z}^2=\mathbb{Z}\langle[F], [B]\rangle$
(use the fibration $*\simeq \mc {\rm P}^2\setminus \mc {\rm P}^1\hookrightarrow \mathcal{F}_3^{2,2}\rightarrow\mathcal{F}_2(\mc {\rm P}^2)$);
the homotopy classes $[F]$ and $[B]$ correspond to the free generators of the second homotopy groups of the fiber and of
the basis respectively, in the fibration (see \cite{B2}) $\mc {\rm P}^1 \simeq (\mc {\rm P}^2\setminus\{*\})
\hookrightarrow\mathcal{F}_2(\mc {\rm P}^2) \rightarrow\mc {\rm P}^2$:
$$ F: (D^2,S^1)\rightarrow(\mathcal{F}_3^{2,2},d_*^0),\,\,\,z\mapsto(d_1^0,d_2^{F(z)}, d_3^{F(z)}), $$
where $d_2^{F(z)}:zX_0-rX_1=0=X_3$ and $d_3^{F(z)}:rX_0+\overline{z}X_1=0=X_3$, and also
$$ B:(D^2,S^1)\rightarrow(\mathcal{F}_3^{2,2},*),\,\, z\mapsto(d_1^{B(z)}, d_2^0, d_3^{B(z)}), $$
where $d_1^{B(z)}:zX_0-rX_3=0=X_1,\,\,d_3^{B(z)}:rX_0+\overline{z}X_3=0=X_1.$
Choosing the lifts $\widetilde{F},\widetilde{B}:(D^2,S^1)\rightarrow(\mathcal{D}_{I^0}^3,\mathcal{F}_2(d_1^0)\times
\mathcal{F}_2(d_2^0)\times\mathcal{F}_2(d_3^0))$:
$$ \widetilde{F}(z)=\big(A_1^0,B_1^0,A_2^{\widetilde{F}(z)},B_2^{\widetilde{F}(z)},A_3^{\widetilde{F}(z)},
B_3^{\widetilde{F}(z)}\big) $$
with
$$ \begin{array}{ll}
A_2^{\widetilde{F}(z)}=[r:z:0:0],             & B_2^{\widetilde{F}(z)}=[r:z:1:0], \\
A_3^{\widetilde{F}(z)}=[\overline{z}:-r:0:0], & B_3^{\widetilde{F}(z)}=[\overline{z}:-r:1:0],
\end{array} $$
respectively
$$\widetilde{B}(z)=\big(A_1^{\widetilde{B}(z)},B_1^{\widetilde{B}(z)},A_2^0,B_2^0,A_3^{\widetilde{B}(z)}
,B_3^{\widetilde{B}(z)}\big)$$
with
$$ \begin{array}{ll}
A_1^{\widetilde{B}(z)}=[r:0:0:z],              & B_1^{\widetilde{B}(z)}=[r:0:1:z]            \\
A_3^{\widetilde{B}(z)}=[\overline{z}:0:0:-r],  & B_3^{\widetilde{B}(z)}=[\overline{z}:0:1:-r],
\end{array} $$
we obtain the equalities $\delta_*([F])=-[b]+[c]$, $ \delta_*([B])=-[a]+[c]$. Therefore we proved that
\begin{corollary}
The fundamental group of the space $\mathcal{D}_I^3$ is infinite cyclic generated by $[\a]$.
\end{corollary}
Using the second fibration of Lemma \ref{lem.4.1}, we find the exact sequence
$$\rightarrow \pi_2(\mc {\rm P}^3)\mathop{\rightarrow}\limits^{\delta_*}\pi_1(\mathcal{D}_{I^0}^3)\rightarrow
\pi_1(\mathcal{D}^3)\rightarrow 1$$
where the generator $\Psi: (D^2,S^1)\rightarrow (\mc {\rm P}^3,I^0),\,z\mapsto[r:0:z:0]$ has the lift
$$\widetilde{\Psi}:(D^2,S^1)\rightarrow\mathcal{D}^3,\,z \mapsto\big(A_1^0,B_1^{\widetilde{\Psi}(z)},
A_2^0,B_2^{\widetilde{\Psi}(z)},A_3^{\widetilde{\Psi}(z)},B_3^{\widetilde{\Psi}(z)}\big),$$
with
$$ \begin{array}{ll}
B_1^{\widetilde{\Psi}(z)}=[r:0:z:1],             & B_2^{\widetilde{\Psi}(z)}=[r:1:z:0],\\
A_3^{\widetilde{\Psi}(z)}=[\overline{z}:0:-r:0], & B_3^{\widetilde{\Psi}(z)}=[r+\overline{z}:0:z-r:0].
\end{array}  $$
Therefore $\delta_*([\Psi])=[\widetilde{\Psi}|S^1]=[\a]+[\b]+2[\g]=4[\a]$, and we proved:
\begin{corollary}
The fundamental group of the space $\mathcal{D}^3$ is cyclic of order four and it is generated by $[\a]$.
\end{corollary}
\begin{proposition}
$$\pi_1(\mathcal{D}_{I}^{3,4})\cong\pi_1(\mathcal{D}_{I}^{3,n})\hspace{2cm}(n\geq4);$$
$$\pi_1(\mathcal{D}^{3,4})\cong\pi_1(\mathcal{D}^{3,n})\hspace{2cm}(n\geq4).$$
\end{proposition}
\begin{proof} This is like in \ref{cor:3.2}.
\end{proof}
\noindent\textit{Proof of Theorem \ref{th:1.4} and of Theorem \ref{th:1.5}}. We show that
$\pi_1(\mathcal{D}_{I}^{3,4})=1$; this implies that $\pi_1(\mathcal{D}^{3,4})=1$. Choose as a generator for
the fundamental group of the space of 3-planes in $ \mc {\rm P}^4 $ containing the fixed point
$ I=[0:0:1:0:0] $ the class of the map
$$\Sigma: (D^2,S^1)\rightarrow ({\rm Gr}^2(\mc {\rm P}^3),\, X_4=0), \,z\mapsto rX_1-zX_4=0. $$
The lift
$$\widetilde{\Sigma}:(D^2,S^1)\rightarrow(\mathcal{D}_{I^0}^{3,4},\mathcal{D}_{I^0}^{3,3}),\,
z \mapsto\big(A_1^{00},B_1^{00},
A_2^{\widetilde{\Sigma}(z)},B_2^{\widetilde{\Sigma}(z)},A_3^{00},B_3^{00}\big),$$
where $ A_1^{00}=[0:0:0:1:0],\dots,B_3^{00}=[1:0:1:0:0] $ are fixed points and
$$ A_2^{\widetilde{\Sigma}(z)}=[0:z:0:0:r],\, B_2^{\widetilde{\Sigma}(z)}=[0:z:1:0:r],  $$
shows that $ \delta_*:\pi_2({\rm Gr}^2(\mc {\rm P}^3))\rightarrow\pi_1(\mathcal{D}_{I^0}^{3,3}) $ is an isomorphism.
\hfill $\square$

\end{document}